\documentclass[12pt,article]{amsart}
\usepackage{amsfonts,amsmath,amssymb}

\newcommand{\C}{{\mathbb C}}

\newcommand\wt[1]{{\widetilde{#1}}}

\newcommand\R{{\mathbb R}}
\newcommand\gS{\Sigma}
\newcommand\gs{\sigma}
\newcommand\Var{\operatorname{Var}}
\renewcommand\mod{\operatorname{mod}}
\newcommand\ch{{\mathcal H}}
\newcommand\gD{\Delta}
\newcommand\cd{{\mathcal D}}
\newcommand\gO{\Omega}
\newcommand\go{\omega}
\newcommand\cb{{\mathcal B}}
\newcommand\gd{\delta}
\newcommand\Cyc{\operatorname{Cyc}}
\newcommand\rank{\operatorname{rank}}

\newtheorem{theorem}{Theorem}

\newtheorem{cor}{Corollary}

\newtheorem{dfn}{Definition}

\begin{document}

\author{
M. M. Malamud and S. M. Malamud}

\title{
On the spectral theory of operator measures
}

\maketitle

{\bf 1. Introduction} Operator measures naturally arise in various questions 
of spectral theory of self-adjoint operators (with spectrum of finite or 
infinite multiplicity), integral representations of operator-valued functions 
of Herglotz and Nevanlinna classes, in the theory of models of 
symmetric operators, etc.    

Throughout the note, $H$ is a separable Hilbert space and $\gS(t)=\gS(t)^*$ 
is a nondecreasing strongly left continuous ($\gS(t-0)=\gS(t))$
operator-function on $\R$ in  $B(H).$ 
In a standart way (see [3], [4]) the function $\gS(t)$ determines an 
operator measure $\gS,$ 
defined on the algebra ${\mathcal B}_b(\R)$ of bounded Borel subsets of  $\R.$ 

The theory of orthogonal measures (resolutions of the identity) is known 
in detail.
In this note, we consider several questions of the theory of 
nonorthogonal operator measures. 
An essential role in our considerations is played by the 
Beresaskii - Gelfand - Kostyuchenko (BGK) theorem on the differentiati on 
of an operator measure [2,3,6]. 

We obtain an inner description of the space 
$L_2(\gS,H).$ This problem was posed by M. G. Krein [9] and (in the special 
case   
dim$H<\infty)$ solved by I. S. Kac 
[1], [7], [8]. 

Further, we 
construct a theory of Hellinger spectral types for 
a nonorthogonal operator measure. We establish the existence 
of subspaces realizing Hellinger spectral types and in particular 
the existence of vectors of maximal type.

Some facts are new even for orthogonal measures and for a finite- 
dimensional space $H$.
We show how the spectral Hellinger types of an operator $A=A^*$ can be found 
via a cyclic subspace $L.$
It turnes out that the set of the vectors of maximal type, lying in $L$ is 
an everywhere dense set of type $G_\gd$ and of second category. 

Moreover, we establish an analog of the Jordan Theorem for Operator 
measures-charges. 
For simplicity,we state all results for measures on the line, even 
though they remain valid for measures defined on subsets of $\R^n.$ 

{\bf 2. The space $L_2(\gS,H).$} 
Following [3], we recall the definition of the space $L_2(\gS,H).$
Let $C_{00}(H)$ be the set of all strongly continuous compactly supported 
vector-functions $f$  
ranging in finite-dimensional 
subspaces of $H$ (the subspace depends on $f).$
Further, for $f,g\in C_{00}(H)$ we introduce the inner product 
$(f,g)_{L_2(\gS,H)}=
\int_\R(d\gS(t)f(t),g(t))_H$.
(The intergal is understood as the limit of Riemann sums).
Faktorizing $C_{00}(H)$  by the lineal $L_0=\{f:\ (f,f)_{L_2(\gS,H)}=0\}$ 
and completing it, we arrive at the Hilbert space $L_2(\gS,H).$

Let  ${\frak S}_2(H)$ be the ideal of Hilbert-Schmidt operators in $H$. 
Next, let $T\in {\frak S}_2(H)$ satisfy ker$T=$ker$T^*=\{0\}$.
By ${\mathcal D}(T^{-1}$ we denote the domain of
$T^{-1}.$
Finally, let $\rho$ be a scalar measure equivalent to  $\gS\ (\gS\sim \rho)$.

By the BGK theorem, the operator measure $\gS_T(\gD):= 
T^*\gS(\gD)T$ is differentiable in the weak sence with respect to $\rho$; 
its density $\Psi(t):=d\gS_T/d\rho(\ge 0)$ exists $\rho$-a.e.  
and ranges in the set ${\mathfrak S}_1(H)$ 
of trace class  operators in $H$. 

Following [5], we can show that the derivative  
$\Psi$ exists $\rho$-a.e. in the norm of ${\mathfrak S}_1(H)$.
(In [5], this was shown for an orthogonal measure $\gS=E$).

Let ${\wt{\frak H}}_t$ be the completion of $\cd(T^{-1})$ with 
respect to the seminorm
$\|h\|_{\wt{\frak H}_t}=\|\Psi(t)^{1/2}T^{-1}f\|.$
By  ${\frak H}_t$ we denote the corresponding quotient space. 
     \begin{theorem}\label{th4.1} Let $T\in{\frak S}_2(H)$ with  
$\ker T=\ker T^*=\{0\},$  nd let $\rho\sim\gS$ be a scalar measure. 
Then the space  $L_2(\gS,H)$
isometrically coincides with the direct intergal 
of the spaces ${\frak H}_t$ by the measure $\rho(t):$ 
      \begin{equation*}
\L_2(\gS,H)=\int_\R\oplus {\frak H}_t d\rho(t)=:{\frak H}.
       \end{equation*}
Moreover, the identity 
    \begin{equation}\label{aga}
\|f\|^2_{L_2(\gS,H)}=\int_\R\|\Psi(t)^{1/2}T^{-1}f(t)\|^2d\rho(t)
     \end{equation}
holds for the dense set of vector-functions $f\in L_2(\gS,H)$   
with values in $H_+=\cd(T^{-1}).$ 
In particular, the space $\frak H$ does not depend on the choice of $T.$ 
   \end{theorem}
If $\gS(t)=E(t)$ is a resolution of identity in $H$ and $f(t)=E(\gD)h$ 
with  $h\in H_+$ and $\gD\in \cb_b(\R)$, then identity \eqref{aga} 
acquires the form 
$$
(E(\gD)h,h)\ (=\|f(t)\|^2_{L_2(\gS,H)})=
\int_\gD\|\Psi(t)^{1/2}T^{-1}h\|^2_Hd\rho(t).
$$
This identity is equivalent to the direct integral form of the 
BGK theorem for the orthogonal measure $E.$ 

For dim$H<\infty$, Theorem 1 is equivalent to the Kac theorem [8],
but our proof is much simplier than any known one. 
Note that there is a essential difference between the cases dim$H<\infty$ 
and dim$H=\infty.$ 

If for dim$H<\infty$ the space $L_2(\gS,H)$ 
can be identified with a some space of $\rho$-measurabel vector functions  
ranging in $H,$ this fails to be true for dim$H=\infty$ even in 
the simplest cases. 

For example, let $\gS_0\ge 0$ be a compact operator in  $H$.  
We set $\gS(t)=0$ for $t\le t_0$ and $\gS(t)=\gS_0$ for $t>t_0.$ 
Then $L_2(\gS,H)=H\_$ is the completion of $H$ 
with respect to the negative norm $\|f\|_{-}=\|\gS_0^{1/2}f\|.$

Furthermore, if the variation of $\gS$ is unbounded, then not all continuous 
compactly supported functions belong to $L_2(\gS,H)$. 

{\bf 3. The multiplicity function of a measure $\gS.$}
     \begin{dfn}\label{dfn1} Let $\rho\sim \gS$, and let 
$\{e_i\}_{i=1}^\infty$ 
be an orthonormal basis of $H.$ Further, let 
$\gs_{ij}(t):=(\gS(t)e_i,e_j),$ $\psi_{ij}(t):=d\gs_{ij}(t)/d\rho$ and
$\Psi_n(t):=(\psi_{ij}(t))_{i,j=1}^n.$ 

The multiplivity function $N_\gS$ and the total multiplicity $m(\gS)$
of an operator measure $\gS,$ are specified by 
    \begin{equation}\label{kratn}
N_\gS(t):=\sup_{n\ge 1}\rank
\Psi_n(t), \ 
\quad m(\gS):= vraisup N_{\gS}(t) \quad   (\mod \rho).
    \end{equation}
    \end{dfn}
The multiplicity function $N_\gS$ is defined $\rho$-a.e. , 
and one can be shown, that it is independent of the choice of 
a basis $\{e_i\}_1^\infty.$ 
    \begin{dfn} \label{dfn2} 
(a) Let $\gS_1$ and $\gS_2$ be operator measures on $\R.$

A measure  $\gS_1$ is said to be subordinated to $\gS_2\ (\gS_1\prec \gS_2)$ 
if $\gS_1$  is absolutely continuous with respect to $\gS_2,$
that is. $\gS_1(\gd)=0$ whenever $\gS_2(\gd)=0.$

(b) We say, that $\gS_1$ is spectrally subordinated to 
$\gS_2$ $(\gS_1 \prec\prec\gS_2)$ if $\gS_1\prec\gS_2$  
and $N_{\gS_1}(t)\le N_{\gS_2}(t)\ (\mod \gS_2).$ 
The measures $\gS_1$ and $\gS_2$ are said to be spectrally equivalent if 
$\gS_1\prec\prec \gS_2$ and  $\gS_2\prec\prec \gS_1.$
      \end{dfn}

Let  $A$ be a selfadjoint operator in $H,$ 
$E(t):=E_A(t)$ its resolution of identity and $L$ a subspace of $H$.
By $H_L$ we denote the minimal $A$-invariant subspace containing 
$L:$\ $H_L=\text{span}\{E(\gd)L: \gd \in {\mathcal B}(\R)\}.$ 
The subspace $L$ is said to be cyclic
$(L\in \Cyc(A))$ if $H_L=H.$ 
If $L=\{\lambda g:\ =\lambda \in \C\},$ then we write $H_g:=H_L.$ 

The following theorem can be proved with the help of Theorem 1.
       \begin{theorem}\label{th4.4} 
Let $\gS$ be a generalized resolution of the identity in $\ch,$ so that  
$\gS(-\infty)=0$, and $\gS(+\infty)=I_{H_1}.$ 
Let $A$ be a selfadjoint operator in $H$ and $E(t)$ the corresponding 
resolution 
of identity. Then the following assertions are true: 

(a) $N_E(t)$ in Definition \ref{dfn1} coincides with the classical 
multiplicity function of $E$ in the sence of 
[4], [10].

(b) $\gS\prec\prec E$ if and only if there exists a Hilbert space 
${\wt H}\supset \ch$ and a unitary operator $U:\, H\to {\wt H}$ such that 
 $\gS(t)=P_{\ch} UE(t)U^*\lceil \ch$ 
(here $P_{\ch}$ is the orthoprojection in $\wt H$ on $\ch$ and
$U^*\lceil\ch$ is restriction of $U^*$ to $\ch).$ 

(c) If  $U^*\ch\in \Cyc A$. 
then $\gS$ is spectrally equivalent to $E.$ 

Conversely, if  $\gS$ is spectrally euivalent to  $E$ and $N_E(t)$ is
$E$-a.e. finite (for example, if $m(E)<\infty$), then $U^*\ch\in \Cyc A.$

In particular (for $H=\wt H$ and $U=I$), the measure  $\gS$ and its 
minimal orthogonal dilation $E$ are spectrally equivalent. 

(d) The resolution of identity $E_Q$ of the operator 
$Q:\ f\to xf$ in $L_2(\gS,H)$ is a minimal orthogonal dilation of $\gS.$ 
      \end{theorem}

Theorem \ref{th4.4} supplemets the known Najmark theorem [1], 
providing an answer to the question as to which resolution of identity 
can be a dilation of $\gS.$ 
      \begin{cor} The multiplication operators $Q_i:\ f\to xf$ 
in the spaces
$L_2(\gS_i,H_i)\ (i=1,2)$ are unitary equivalent iff $\gS_1$ 
and $\gS_2$ are spectrally equivalent. 
   \end{cor}
{\bf 4. Elements of maximal type.}
Every operator measure $\gS$ in $H$ generates a family of $\gs$-finite 
scalar measures 
$\mu_f\ (\mu_f(\delta):=$\ $(\gS(\delta)f,f))$ 
defined on the algebra 
${\mathcal B}_b(\R).$ 
It is clear that $\mu_f\prec \gS$ for all $f\in H.$ 
It is known ([1],[4]) that any orthogonal measure $E$ in  $H$ posesses 
an element $f$ of maximal type, that is, an element such that $\mu_f\sim E.$ 

This remains valid for nonorthogonal measures. 
Moreover, the following stronger conjecture is valid. 
We note that this conjecture is new even for orthogonal measures, as well
as for the case in which $dim  H<\infty.$ 

       \begin{theorem}\label{main} 
Let  $\gO_\gS := \{f\in H: \gS\sim \mu_f\}$
be the set of all vectors of maximal type for an operator measure $\gS$ 
in $H.$ 
Then: 

\noindent (a) $H\setminus\gO_\gS$ is an $F_\gs$-set of first category in $H;$

\noindent (b) $\omega(\Omega_\gS)=1$ for any Gaussian measure $\omega$ i $H.$ 
    \end{theorem}
       \begin{cor}\label{thcyc} Let  $A$ be a selfadjoint operator in  
$H, let E(t)$ be its resolution of identity, and let $L\in\Cyc A.$ 
Then:

(a) $L\setminus\gO_E$ is an $F_\gs$-set of first category in $L;$

(b) $\go(\gO_E\cap L)=1$ for any Gaussian measure $\omega$ in $L.$ 
      \end{cor}
{\bf 5. Hellinger types.}
The class of all Borel measures, equivalent to a measure $\mu$ 
is called the type of the measure $\mu$ and is denoted by $[\mu]$ (see [4]). 
Let  $A=A^*,\ E:=E_A,\ g\in H$ and 
$\mu_g: \delta \to \mu_g(\delta):=(E(\delta)g,g),\  
\delta \in {\mathcal B}(\R).$  
The type $[g]$ of an element $g$ (with respect to $E)$ is the type 
of the measure 
$\mu_g,\  [g]=[\mu_g].$ 

Consider an orthogonal decomposition of the form  
$H=\oplus_{i=1}^m H_{g_i},\ (m\le \infty).$ 
 If the types of the elements $g_i$ do not increase,
$[g_{i+1}]\prec [g_i],$ then 
their number $m(\le \infty)$ and types
are uniquely determined and are referred to as the Hellinger types 
of the measure $E.$
They form (see [4]) a complete set of unitary invariants of the operator $A.$

Let  $g_1\in \gO_E := \{g\in H:\ \mu_g \sim E\}$.
Then $\mu_g\prec \mu_{g_1}:=\mu$ and for all $g\in H$ the type $[g]$ 
is uniquely determined by the support  
$\Gamma(g):=\{t\in \R:\ d\mu_g/d\mu>0\}$ of the measure 
$\mu_g$ with respect to $\mu$. 
Therefore, the Hellinger types are uiquely determined (mod $E)$ by 
their supports  $\Gamma_i(E) := \Gamma(g_i), i\le m.$

The sets $\Gamma_i(E)$ themselves are  determined (see [4]) by 
the multiplicity function (and the measure $\mu$):  
$\Gamma_i(E)=\{t\in \R:\, N_E(t)\ge i\}.$ 

The existence of a multiplicity function $N_\gS$ of the form 
\eqref{kratn} allows  us to introduce the $i$-th Hellinger type for 
a nonorthogonal measure $\gS$ 
as the type of the scalar measure $d\mu_i:=\chi_id\rho$ with $\rho\sim\gS$ and
$\chi_i$ the indicator of the set 
   \begin{equation}
\Gamma_i(\gS)=\{t\in \R:\, N_\gS(t)\ge i\}, \qquad i\in\{1,\dots, m(\gS)\}. 
     \end{equation}
We refer to $\Gamma_i(\gS)$ as the support of the 
$i$-th Hellinger type of the measure $\gS.$ 
It is clear that 
$\Gamma_i(\gS) \supset \Gamma_{i+1}(\gS).$  Let $i_0$ 
be the number of $\Gamma_i,$ equivalent to $\Gamma_1(\gS)$ $(\mod \rho),$ 
that is, $\rho(\Gamma_1(\gS)\setminus \Gamma_i(\gS))=0.$

If $g_1\in \gO_\gS$ (by Theorem 3 $\gO_\gS\not =\emptyset$),
then $\mu_g\prec \mu_{g_1}=:\mu\sim \gS$ 
with $\mu_g(\gd):=(\gS(\gd)g,g).$
Therefore, the set  $\Gamma(g):=\{t:\ d\mu_g/d\mu>0\}$ is a 
(nontopological) support of the measure $\mu_g$ (that is
$\mu_g(\R\setminus \Gamma(g))=0).$ 

It turns out that although $\gO_\gS\not =\emptyset$,
elements of ``junior'' types  
( that is, vectors $g\in H\setminus \gO_\gS,$ such that
$\Gamma(g)=\Gamma_i(\gS)$ for some $i>i_0)$  may fail to exist.

For example,let $\gS(t)=\sum_{i<t}P_i$ be a $2\times 2$ discrete measure 
with jumps $P_1=\begin{pmatrix} 1&0\\0&0\end{pmatrix},\ 
P_2=\begin{pmatrix} 0&0\\
0&1\end{pmatrix}, \ 
P_3=\begin{pmatrix} 1&0\\0&1\end{pmatrix}  $
at the points $1,2,3.$ 
One can readily see that
$N_\gS(1)=N_\gS(2)=1$ and $N_\gS(3)=2.$ 
Therefore, $\Gamma_1(\gS)=\{1,2,3\}$ and $\Gamma_2(\gS)=\{3\}.$
If $h=(h_1,h_2),$ 
then $\Gamma(h)=\{1,2,3\}=\Gamma_1(\gS)$ 
for $h_1h_2\not=0.$ 
Further, 
$\Gamma(h)=\{2,3\}$ if $h_1=0$  
and $\Gamma(h)=\{1,3\}$ if $h_2=0.$  
Therefore, $\Gamma(h)\not=\{3\}=\Gamma_2(\gS)$ for any $h.$

The desire to realize the ``junior''Hellinger types with the help 
of some subspaces forced us to introduce the following definition.
     \begin{dfn} A subspace $L=L_k\ (\dim L_k=k)$   
is called a $k$-th Hellinger subspace for an operator measure 
$\gS$ in $H$ if $\Gamma_i(P_L \gS \lceil L)=\Gamma_i(\gS)$ 
for all $i\le k$. (Here $P_L$ is the orthoprojection on $L$). 
    \end{dfn}
In particular, a vector of maximal type generates a one-dimensional 
(first) Hellinger subspace. 
       \begin{theorem}\label{main2} 
Under the assumptions of  Theorem \ref{main}, for each $h\in \gO_\gS$ 
there exists a chain of Hellinger subspaces (Hellinger chain):
          \begin{equation}\label{hell}
\{\lambda h\}
=:H_1\subset H_2\subset \ldots \subset H_k \subset \ldots \subset H_m,
\quad \dim H_k=k.
          \end{equation} 
   \end{theorem}
Let  $H_1\subset\ldots\subset H_m$ be a chain of subspaces of the form  
\eqref{hell}, and let $\{e_i\}_1^m$ be a basis in ¢ $H_m,$ such that 
$H_k=$span$\{e_i\}_1^k,\ k\in\{1,\ldots, m\}.$ We set
$\gs_{ij}(t):=(\gS(t)e_i,e_j),$ $\psi_{ij}(t):=d\gs_{ij}(t)/d\rho$, and 
$\Psi_k(t):=(\psi_{ij}(t))_{i,j=1}^k$, where $\rho \sim \gS$. 
The chain $H_1\subset\ldots\subset H_m$ is a Hellinger chain iff 
   \begin{equation} 
\Gamma_k(\gS) = \{t\in \R:\ \det \Psi_k(t) \not=0 \}\ (\mod\gS),
\quad k\in \{1,\dots, m\}.
    \end{equation}

Thus, Theorem 4 amounts to saying that there exists an orthonormal system 
 $\{e_i\}_1^m$ in $H$ such that the $k$-th Hellinger type is realized 
by the measure
$(\land^k\Psi (t)\varphi _k,\varphi _k)d\mu,$ where $\land^k\Psi(t)$
stands for the $k$-th exterior power of  $\Psi(t):=d\gS_T/d\mu$ and
$\varphi_k:=e_1\land \ldots \land e_k$ is a $k$-vector, 
$\varphi_k\in \land^k(H).$   
       \begin{cor}\label{cor2} Let  $A$ be a selfadjoint operator in $H$, 
$E(t)$ its resolution of identity, $L\in\Cyc A$, and 
$h\in \gO_E\cap L.$. 
Then there exists a chain of Hellinger subspaces in $L$ of the form 
\eqref{hell}, $H_k\subset L,\ k\in \{1,\ldots,m\}.$ 

If the multiplicity  $m=m(E)$ of the measure $E$ is finite, then 
$H_m\in \Cyc A.$ 
      \end{cor}
If the vectors $\{e_i\}_{i=1}^m$ 
are pairwice spectrally orthogonal with respect to $E$ and 
$H_k:=$span$\{e_i\}_1^k,$ then the chain 
$H_1\subset\ldots\subset H_m$ is a Hellinger chain if and only if 
$\Gamma(e_i)=\Gamma_i(E),\ i\le m.$ 

If $L\in \Cyc A,$ then in general there is no system of spectrally 
orthogonal vectors realizing 
Helliger types.  
But, according to Corollary 3, these types are realized by means of subspaces 
$H_k\subset L$ or, equivalently, the $k$-vectors 
$\varphi_k=e_1\land \ldots \land e_k \in \land ^k(L).$   

{\bf 6. An analog of the Jordan Theorem.} 
The results of this subsection belong to the second author.
    \begin{dfn} One says that an operator-function 
$\gS=\gS^*:\ \cb(\R)\to B(H)$
(an operator measure-charge) is of weakly bounded variation on $\R,$ if 
 $\mu_{f,g}:\ \gd\to(\gS(\gd)f,g)$ 
is a finite charge on the Borel  $\sigma$-algebra $\cb(\R)$ 
for any $f,g\in H.$ 
   \end{dfn}
      \begin{theorem} An operator measure-charge $\gS$ of weakly bounded variation on $\R$ can be be expressed as the difference of two finite nonnegative 
operator measures $\gS=\gS_1-\gS_2$ 
if and only if 
$$
\Var_\R\|T^*\gS(\gD)T\|_1:=\sup_{\pi}\sum_i\|T^*\gS(\gD_i)T\|_1
=:c(T)<\infty\qquad 
$$
for any $T\in{\frak S}_2(H),$  where $\|\cdot\|_1$ is the trace norm 
and the supremum is taken over all partitions $\pi=\{t_j\}_{-\infty}^{\infty}$
of \ $\R,$  \  $\gD_i=[t_i,t_{i+1}).$ 
       \end{theorem}
The proof is based on some facts from the theory of $C^*-$algebras  
and completely bounded maps [12].

Let  $x_i^{(n)}=(x_i^{(n)})^*$ be the Clifford $2^n\times 2^n$-matrices, 
$x_i^{(n)}x_j^{(n)}+x_j^{(n)}x_i^{(n)}=2\gd_{ij}I,\  i,j \in \{1,\ldots, n\}.$
Consider the operator measure-charge
$$
\gS(\gD)=\oplus_1^\infty \gS_n(\gD),\qquad 
\gS_n(\gD)=1/\sqrt{2n}\sum_{1/k\in \gD,\ k\le n}x_k^{(n)}.
$$ 
Clearly, 
$\gS(\gD)$ is a discrete measure with support 
$\{0\}\cup\{1/k\}_{k=1}^\infty.$ 
It is of weakly bounded variation but can not be expressed as the difference 
of two nonnegative operator measures. 

The idea of using Clifford matrices for the constraction of the measure
$\gS$ has been borrowed from [10], where it was used for another purpose.

We express our gratitude to Yu. M. Berezanskii and M. Sh. Birman 
for useful discussions, which improved the paper. 

\begin{center}
{\sc References }
\end{center}

{\bf 1.} {\it N. I. Achiezer and I. M. Glazman,} 
The theory of linear operators in Hilbert space, 
Nauka, M., 1966. 
{\bf 2.} {\it Yu. M. Berezanskii,} Dokl. Acad. Sci. USSR, {\bf 108}, 
No. 3, 379-382 (1956).
{\bf 3.} {\it Yu. M. Berezanskii,} Expansions in eigenfunctions of selfadjoint operators, Naukova Dumka, K. 1965.
{\bf 4.} {\it M. Sh. Birman and M. Z. Solomyak,} Spectral theory of selfadjoint operators in Hilbert space, Len. Univ., 1980.
{\bf 5.} {\it  M. Sh. Birman and C. B. Entina,} Izv. AN. USSR, {\bf 31}, 401-430 (1967).
{\bf 6.} {\it I. M. Gelfand and A. G. Kostyuchenko,} Dokl. Acad. Sci USSR, {\bf 103}, No. 3, 349-352 (1955).
{\bf 7.} {\it N.Dunford and J. T. Schwarz,} Linear operators, vol. 2
{\bf 8.} {\it I. S. Kac,} Zap. Chark. Mat. soc. (4), 22, 95-113 (1950).
95-113 (1950).
{\bf 9.} {\it M. G. Krein,} Trans. math. inst. AN USSR, {\bf 10}, 83-106 (1948).
{\bf 10.} {\it V. Paulsen,} J. Funct. Anal. {\bf 109}, 113-129 (1992).
{\bf 11.} {\it A.I. Plesner,} Spectral Theory of operators, Nauka, M.1965.
{\bf 12.}  {\it G. Pisier,} Similarity problems and completely
bounded maps, Springer, 1996. 
\end{document}